\documentclass{amsart}

\usepackage{epsfig,amssymb}

\newcommand\ns[1]{\mbox{\small $#1$}}

\newtheorem{lem}{Lemma}

\begin{document}

\title{\bf 
On the polynomial moment problem.
}
\author {F. Pakovich}
\address{Department of Mathematics, Weizmann Institute of Science,
Rehovot 76100, Israel}
\email{pakovich@wisdom.weizmann.ac.il}


\date{\today}


\maketitle

\def\d{{\rm d}}
\def\I{{\rm I}}

\def\C{{\mathbb C}}
\def\N{{\mathbb N}}
\def\P{{\mathbb P}}
\def\D{{\mathbb D}}
\def\Z{{\mathbb Z}}
\def\phi{{\varphi }}
\def\d{{\rm d\,}}
\def\deg{{\rm deg\,}}
\def\Det{{\rm Det}}\def\dim{{\rm dim\,}}
\def\Ker{{\rm Ker\,}}
\def\Gal{{\rm Gal\,}}
\def\St{{\rm St\,}}
\def\Sym{{\rm Sym\,}}
\def\Mon{{\rm Mon\,}}

\vskip 0.2cm

\section{Introduction.}

In this paper we treat the following "polynomial moment problem":
{\it for a complex polynomial $P(z)$ and distinct $a,b \in \C$
such that $P(a)=P(b)$ to describe polynomials $q(z)$
such that
$$\int^b_a P^i(z)q(z)\d z=0 \eqno(*)$$
for all integer non-negative $i.$
}

The polynomial moment problem was
proposed in the series of papers of M. Briskin, J.-P. Francoise and
Y. Yomdin \cite{bfy1}-\cite{bfy5} as an infinitesimal version of
the center problem for the polynomial Abel equation in the complex domain
in the frame
of a programme concerning the classical Poincare center-focus
problem for the polynomial vector field on the plane.
It was suggested that the following "composition condition" 
imposed on $P(z)$ and $Q(z)=\int q(z)\d z$ is necessary
and sufficient for the pair $P(z),$ $q(z)$ to satisfy (*):
{\it there exist polynomials $\tilde P(z),
\tilde Q(z), W(z)$ such that} $$P(z)=\tilde P(W(z)), \ \ \
Q(z)=\tilde Q(W(z)), \ \ \ {\it and} \ \ \ W(a)=W(b).\eqno(**)$$

It is easy to see that the composition condition is sufficient:
since after the change of variable $z\rightarrow W(z)$ 
the way of integration becomes closed, the sufficientness
follows from the Cauchy theorem.
The necessity of the composition condition in the case when $a,b$ are not
critical points of $P(z)$
was proved by C. Christopher
in \cite{c} (see also the paper of N. Roytvarf \cite{ro}
for a similar result) and in
some other special cases by M. Briskin, J.-P. Francoise and
Y. Yomdin in the papers cited above.
Nevertheless, in general the composition conjecture fails to be true
as it was shown by the author in \cite{pa}. 

In this paper we give a solution of the polynomial moment problem
in the case when $P(z)$ is indecomposable that is when $P(z)$
can not be represented as a composition
$P(z)=P_1(P_2(z))$ with non-linear polynomials $P_1(z),$ 
$P_2(z).$ We show that in this case the
composition conjecture is true without any restrictions on points $a,b.$

\pagebreak

\noindent{\bf Theorem 1.} {\it 
Let $P(z),q(z)$ be complex polynomials and let $a,b$ be distinct complex
numbers such that $P(a)=P(b)$ and
$$\int^b_a
P^i(z)q(z)\d z=0$$
for $i\geq 0.$ Suppose that $P(z)$ is indecomposable. Then there exists a
polynomial $\tilde Q(z)$ such that
$Q(z)=\int q(z)\d z=\tilde Q(P(z))$.}

\vskip 0.2cm

We also examine the following condition which is stronger than (*):
$$\int^b_a
P^i(z)Q^{j}(z)Q^{\prime}(z)\d z=0$$ for $i\geq 0,$
$j\geq 0.$
If $\gamma$ is a curve which is the image of the segment $[a,b]$
in $\C^2$ under the map $z\rightarrow (P(z),Q(z))$
then this condition is equivalent
to the condition that $\int_{\gamma}\omega =0$ for all global
holomorphic 1-forms $\omega$ in $\C^2$ ("the moment condition").
For an oriented simple closed curve $\delta$ of class $C^2$ in $\C^2$ 
the moment condition is necessary
and sufficient to be a boundary of a bounded analytic variety
$\Sigma$ in $\C^2;$ it is a special case of
the result of R. Harwey and B. Lawson \cite{hl}. 
The case when $\delta$ is an image
of $S^1$ under the map $z \rightarrow (f(z), g(z)),$ where $f(z),g(z)$
are functions analytic in an annulus containing $S^1$ was investigated
earlier by J. Wermer \cite{w1}:
in this case the moment condition is equivalent to the condition that
there exists a finite Riemann surface $\Sigma$ with border $S^1$
such that $f(z),$ $g(z)$ have an 
analytic extension to $\Sigma.$

Unlike to condition (*) the more restrictive moment
condition imposed on polynomials $P(z),$ $Q(z)$ 
turns out to be equivalent to composition condition (**). 
We show that actually even weaker 
condition is needed.

\vskip 0.2cm

\noindent{\bf Theorem 2.} {\it 
Let $P(z),Q(z)$ be complex polynomials and let $a,b$ be distinct complex
numbers such that $P(a)=P(b)$ and 
$$\int^b_a P^i(z)Q^j(z)Q^{\prime}(z)\d z=0$$
for $0\leq i\leq \infty,$ $0\leq j \leq d_a+d_b-2,$
where $d_a$ (resp.
$d_b$) is the multiplicity of the point $a$ (resp. $b$) with respect to
$P(z).$ Then there exist polynomials $\tilde P(z),$
$\tilde Q(z),$ $W(z)$ such that $P(z)=\tilde P(W(z)),$
$Q(z)=\tilde Q(W(z)),$ and $W(a)=W(b).$}

\vskip 0.2cm

Note that if $a,b$ are not critical points of $P(z)$
that is if $d_a=d_b=1$ then conditions of the theorem reduce to 
condition (*) and therefore theorem 2
includes as a particular case the result of C. Christopher.

\section{Proofs.}

\subsection{Lemmata about branches of $Q(P^{-1}(z))$.}

\vskip 0.2cm

\noindent Let $P(z)$ and $Q(z)$ be rational functions and
let $U\subset \C$ be a domain in which there exists a single-valued branch
$p^{-1}(z)$ of the algebraic function $P^{-1}(z).$
Denote by $Q(P^{-1}(z))$
the complete algebraic function obtained by the analytic
continuation of the functional element $\{U,Q(p^{-1}(z))\}.$
Since the monodromy group $G(P^{-1})$ of the algebraic function
$P^{-1}(z)$ is transitive this definition does not depend
of the choice of $p^{-1}(z).$ Denote by $d (Q(P^{-1}(z)))$ the degree of
the algebraic function $Q(P^{-1}(z))$ that is the number of its
branches.

\begin{lem}

Let $P(z),Q(z)$ be rational functions.
Then $$d (Q(P^{-1}(z)))=\deg P(z)/[\C(z):\C(P,Q)].$$

\end{lem}

\noindent{\it Proof.} Since any algebraic relation over $\C$ between
$Q(p^{-1}(z))$ and $z$ supplies
an algebraic relation between $Q(z)$ and $P(z)$ and vice versa
we see that $d (Q(P^{-1}(z)))=[\C(P,Q):\C(P)].$ As
$[\C(P,Q):\C(P)]=[\C(z):\C(P)]/[\C(z):\C(P,Q)]$ the lemma follows
now from the observation that $[\C(z):\C(P)]=\deg P(z).$

\vskip 0.2cm

Recall that by L\"{u}roth theorem 
each field $k$ such that $\C \subset k \subset \C(z)$ and
$k\neq \C$ is of the form $k=\C(R),$ $R\in \C(z)\setminus \C.$ Therefore,
the field $\C(P,Q)$ is a proper subfield of $\C(z)$ if and only if
$P(z)=\tilde P(W(z)),$ $Q(z)=\tilde Q(W(z))$
for some rational functions $\tilde P(z),$ $\tilde Q(z),$
$W(z)$ with $\deg W(z) > 1;$ in this case we say that $P(z)$ and
$Q(z)$ have a common right divisor in the composition algebra.
The lemma 1 implies the following explicit criterion which
essentially due to Ritt \cite{ri} (cf. also \cite{c}, \cite{ro}).

\vskip 0.2cm

\noindent{\bf Corollary 1.} {\it Let $P(z),Q(z)$ be rational functions.
Then $P(z)$ and $Q(z)$ have a common right divisor in
the composition algebra if and only if
$$Q(p^{-1}(z))=Q(\tilde p^{-1}(z))\eqno(1)$$
for two different branches $p^{-1}(z),$ $\tilde p^{-1}(z)$
of $P^{-1}(z).$}

\vskip 0.2cm

\noindent{\it Proof.} Indeed, by lemma 1 the field $\C(P,Q)$ 
is a proper subfield of $\C(z)$ if and only if $d (Q(P^{-1}(z)))<\deg
P(z).$ On the other hand, the last inequality is clearly equivalent to
condition (1).

\begin{lem}

Let $P(z),$ $Q(z)$ be rational functions, $\deg P(z)=n.$ 
Suppose that there exist
$a_i\in \C,$ $1\leq i \leq n,$ not all equal between themselves such that
$$\sum_{i=1}^na_iQ(p^{-1}_i(z))=0. \eqno(2)$$
If, in addition, the group $G(P^{-1})$ is doubly transitive 
then $Q(z)=\tilde Q(P(z))$ for some rational function $\tilde Q(z).$ 

\end{lem}

\noindent{\it Proof.}
Let $G\subset S_n$ be a permutation
group and let $\rho_G: G\rightarrow GL(\C^n)$ be
the permutation representation of $G$  
that is $\rho_G(g),$ $g\in G$ is the linear map which sends 
a vector $\vec a=(a_1,a_2, ... , a_n)$ to the vector
$\vec{a_g}=(a_{g(1)},a_{g(2)}, ... , a_{g(n)}).$
It is well known
(see e.g. \cite{wi}, Th. 29.9) that $G$ is doubly transitive if and 
only if $\rho_G$ is
the sum of the identical representation and an absolutely irreducible
representation. Clearly, the one-dimensional $\rho_G$-invariant subspace
$E\subset \C^n$
corresponding to the identity representation is generated by the vector
$(1,1, ... ,1).$ Therefore, since the Hermitian inner product
$(\vec{a},\vec{b})=a_1 \bar b_1+ a_2 \bar b_2 + ... + a_n \bar
b_n$ is invariant with respect to $\rho_G,$ the group $G$ is doubly
transitive if and only if the subspace $E$ and its orthogonal complement
$E^{\perp}$ are the only $\rho_G$-invariant subspaces of $\C^n.$

Suppose that (2) holds. In this case  
also $$\sum_{i=1}^na_iQ(p^{-1}_{\sigma(i)}(z))=0 \eqno(3)$$ for all
$\sigma\in G(P^{-1})$ by the analytic continuation. To prove the lemma 
it is enough to show that
$Q(p^{-1}_{i}(z))=Q(p^{-1}_{j}(z))$ for all $i,j,$ $1\leq
i,j \leq n;$ then by lemma 1
$[\C(z):\C(P,Q)]= \deg P(z)=[\C(z):\C(P)]$ and therefore
$Q(z)=\tilde Q(P(z))$ for some rational function
$\tilde Q(z)$. Assume the converse i.e. that there exists 
$z_0\in U$ such that not all $Q(p^{-1}_{i}(z_0)),$
$1\leq i \leq n,$ are equal between themselves. 
Without loss of generality we can suppose that all $Q(p^{-1}_{i}(z_0)),$
$1\leq i \leq n,$ are finite.
Consider the subspace
$V\subset \C^n$ generated by the vectors $\vec v_{\sigma},$
$\sigma\in G(P^{-1}),$ where $\vec v_{\sigma}=(Q(p^{-1}_{\sigma(1)}(z_0)),
Q(p^{-1}_{\sigma(2)}(z_0)), ...
, Q(p^{-1}_{\sigma(n)}(z_0)).$ Clearly, $V$ is $\rho_{G(P^{-1})}$-invariant
and $V\neq E.$
Moreover, it follows from (3) that $V$ is
contained
in the orthogonal complement $A^{\perp}$ of the subspace $A\subset \C^n$
generated by the vector $(\bar a_1, \bar a_2, ... , \bar a_n).$ 
Since $A\neq E$ we see that $V$ is a
proper $\rho_G$-invariant subspace of $\C^n$ distinct from $E$ and
$E^{\perp}$ that contradicts the assumption that
the group $G(P^{-1})$ is doubly transitive.

\subsection{Lemma about preimages of domains.}
For a polynomial $P(z)$ denote by $c(P)$ the set of finite critical
values of $P(z)$. 

\begin{lem}

Let $P(z)$ be a polynomial and let $U\subset \C\P^1$ be an unbounded simply 
connected domain such that $c(P)\cap U=\emptyset.$ 
Then $P^{-1}\{U\}$ is conformally
equivalent to the unit disk and $P^{-1}\{\partial U\}$ is connected.

\end{lem}

\noindent{\it Proof.} Indeed, by the Riemann theorem $U$
is conformally equivalent to the unit disk $\D$ whenever $\partial
U$
contains more than one point. It follows from $c(P)\cap U=\emptyset$
that $\partial U$ contains a unique point
if and only if $P(z)$ has a unique finite critical value $c$ and
$\partial U=c;$ in this case there exist linear functions
$\sigma_1,$ $\sigma_2$ such that $\sigma_1(P(\sigma_2(z)))=z^n,$ 
$n\in \N$ and the lemma is obvious.
Therefore, we can suppose that $U\cong \D.$ Since $c(P)\cap U=\emptyset$
the restriction of the map
$P(z)\,:\, \C\P^1\rightarrow \C\P^1$ on $
P^{-1}\{U\}\setminus P^{-1}\{\infty\}$ is a covering map.
As $U \setminus\infty $
is conformally equivalent to the punctured
unit disc $\D^*$ it follows from covering spaces
theory that $ P^{-1}\{U\}\setminus
P^{-1}\{\infty\}$ is a disjoint union of domains $\cup U_i$
conformally equivalent to $\D^*$ such that all induced maps
$f_i \,:\, \D^* \rightarrow \D^*$ are of the form
$z\rightarrow z^{l_i},$ $l_i\in \N.$
But, as $P^{-1}\{\infty\}=\{\infty\},$ there may be only one such a
domain. Therefore, the preimage $P^{-1}\{U\}$ is conformally
equivalent to the unit disk. In particular, since
$P^{-1}\{\partial U\}=
\partial P^{-1}\{U\}$ we see that $P^{-1}\{\partial U\}$ is
connected.

\subsection{Proof of theorem 2: the case of a regular value.}

In this section we investigate the case
when $t_0=P(a)=P(b)$ is not a
critical value of the polynomial $P(z).$
For a simple closed curve $M\subset \C$ denote by $D_M^+$ (resp. by
$D_M^-$) the domain that is interior (resp. exterior) with respect to $M.$

Let $L\subset \C$ be a simple closed curve
such that $t_0\in L$ and $c(P)\subset D_L^+.$ Denote by $\vec{L}$ the
same curve considered as an embedded into the complex plane oriented
graph. By definition, the graph $\vec{L}$ has one vertex $t_0$ and one
counter-clockwise oriented edge
$l.$ Let $\vec{\Omega}=P^{-1}\{\vec{L}\}$ be an oriented graph which is
the preimage of the graph $\vec{L}$ under the mapping $P(z): \C
\rightarrow \C,$ i.e. vertices of $\vec\Omega$ are preimages of $t_0$ and
oriented edges of $\vec\Omega$ are preimages of $l.$ 
As $L\cap c(P)=\emptyset$ the graph $\vec\Omega$ has $n=\deg P(z)$ vertices
and $n$ edges. Furthermore, by lemma 3
the graph $\vec\Omega=P^{-1}\{\partial D_L^-\}$ is connected.
Therefore, as a point set in $\C$ the graph $\vec\Omega$ is
a simple closed curve. Let $l_j,$ $1\leq j \leq n,$ be oriented
edges of $\vec\Omega$ and let $a_j$ (resp. $b_j$) be the
starting (resp. ending) point of
$l_j.$ We will suppose that edges of $\vec\Omega$
are numerated by such a way that $a_1=a$ and that under a moving
around the domain $P^{-1}\{D_L^-\}$ along its boundary $\vec\Omega$
the edge $l_{i},$ $1\leq i \leq n-1,$ is followed by the edge $l_{i+1}$
(see fig. 1).

Let $U\subset \C$ be a simply connected domain such that
$U\cap c(P)=\emptyset$ and $L\setminus \{t_0\}\subset U.$ By monodromy
theorem, in such a domain there exist $n$ single-valued branches of
$P^{-1}(t).$ Denote by $p^{-1}_j(t),$ $1 \leq j \leq n,$
the single-valued branch of $P^{-1}(t)$ defined in $U$ by the
condition $p_j^{-1}\{l\setminus t_0\}=l_j\setminus\{a_j,b_{j}\};$
such a numeration of branches of $P^{-1}(t)$
means that the analytic continuation of the functional
element $\{U,p_j^{-1}(t)\},$ $1 \leq j \leq n-1,$ along $L$
is the functional element $\{U,p_{j+1}^{-1}(t)\}.$
Let $l_k,$ $k<n,$ be the edge of $\Omega$ such that $b_k=b$ and
let $\Gamma=\{l_1,l_2, ... ,l_k\}$ be the oriented path in the graph
$\Omega$ joining the vertices $a_1=a$ to $b_k=b.$
For $t\in U$ set
$\phi(t)=\sum_{j=1}^k
Q(p^{-1}_j(t)).$ Clearly, $\phi(t)$ is analytic on $U$ and extends to a 
continuous on $U\cup t_0$ function
since $\sum_{j=1}^k Q(a_j)-\sum_{j=1}^k Q(b_j)=Q(a)-Q(b)=0.$

\begin{figure}[t]
\centering
\input{fig1.pstex_t}
\caption{}
\end{figure}

Consider an analytic on $\C\P^1\setminus L$ function
$$I(\lambda)=\oint_L\frac {\phi(t)}{t-\lambda}\d t=
\int_{\Gamma}\frac{Q(z)P^{\prime}(z)\d z}{P(z)-\lambda}.$$
More precisely, the integral above defines two analytic functions:
one of them $I^+(\lambda)$ is analytic in $D_L^+$ and the other
one $I^-(\lambda)$ is analytic in $D_L^-.$
Furthermore, calculating the Taylor expansion at infinity 
we see that condition (*) reduces
to the condition that $I^-(\lambda)\equiv 0$ in $D_L^-.$
Therefore, by a well-known result about integrals of the
Cauchy type (see e.g. \cite{mu}, p. 63) the function $\phi(t)$
is the boundary value on $L$ of the analytic in $D_L^+$ function 
$I^+(\lambda)$. It follows from the uniqueness theorem for boundary values
of analytic functions that the functional element 
$\{U,\phi(t)\}$ can be analytically continued along any curve 
$M\subset D_L^+.$
As $c(P)\subset D_L^+$ this fact implies that 
$\{U,\phi(t)\}$ can be analytically continued along any curve 
$M\subset \C.$
Therefore, by the monodromy theorem, the element $\{U,\phi(t)\}$ 
extends to a single-valued analytic function in the whole complex plane.
In particular,
the analytic continuation of $\{U,\phi(t)\}$ along any closed
curve coincides with 
$\{U,\phi(t)\}.$ On the other hand, by construction the analytic
continuation of $\{U,\phi(t)\}$ along the curve $L$ is
$\{U, \phi_L(t)\},$ where $\phi_L(t)=\sum_{j=2}^{k+1}
Q(p_{j}^{-1}(t)).$ It follows from $\phi(t)=\phi_L(t)$ 
that $Q(p_{1}^{-1}(t))=Q(p_{k+1}^{-1}(t))$
and by corollary 1 we conclude that $P(z)$ and $Q(z)$ have a
common right divisor in the composition algebra.

As the field $\C(P,Q)$ is a proper subfield of $\C(z)$ and $P(z),$ 
$Q(z)$ are polynomials it is easy to prove that $\C(P,Q)=\C(W)$ for some
polynomial $W(z),$ $\deg W(z) > 1.$ It means that 
$P(z)=\tilde P(W(z)),$ $Q(z)=\tilde Q(W(z))$ for 
some polynomials $\tilde P(z),$ $\tilde Q(z)$ such that
$\tilde P(z)$ and $\tilde Q(z)$ have no a
common right divisor in the composition algebra.
Let us show that $W(a)=W(b).$
Since $t_0$ is not a critical value of the polynomial $P(z)=\tilde
P(W(z))$ the chain rule implies that $t_0$ is not a critical value of 
the polynomial $\tilde P(z).$
Therefore, if $W(a)\neq W(b)$ then 
after the change of variable $z\rightarrow W(z)$
in the same way as above we find that $\tilde P(z)=\bar
P(U(z))$, $\tilde Q(z)=\bar
Q(U(z))$ for some polynomials $\bar P(z), \bar Q(z), U(z)$ with
$\deg U(z) >1$ that contradicts the fact that
$\tilde P(z),$ $\tilde Q(z)$ have no a
common right divisor in the composition algebra.
This completes the proof in the case when $z_0$ is not a critical value
of $P(z).$

\begin{figure}[t]
\centering
\input{fig2.pstex_t}
\caption{}
\end{figure}

\subsection{Proof of theorem 2: the case of a critical value.}
Assume now that $t_0=P(a)=P(b)$ is a critical value of
$P(z).$ In this case let $L$ be a
simple closed curve such that $t_0\in L$ and $c(P)\setminus
t_0\subset D_L^+.$ Consider again a graph 
$\vec\Omega=P^{-1}\{\vec{L}\}$. 
Since $P^{-1}\{D_L^-\}$ is still 
conformally equivalent to the unit disk by lemma 3, we see
that the graph $\vec\Omega$ topologically is the boundary of a disc
although it is not a simple closed curve any more.
Let $l_j,$ $1\leq j \leq n,$ be
oriented edges of $\vec\Omega$ and let $a_j$ (resp. $b_j$)
be the starting (resp. the ending)
point of $l_j.$
Let us fix again such a numeration of edges of $\vec\Omega$ that $a_1=a$
and that under a moving around the domain $P^{-1}\{D_L^-\}$ along its
boundary $\vec\Omega$
the edge $l_{i},$ $1\leq i \leq n-1,$ is followed by the edge $l_{i+1}.$
As above denote by $U$ a domain in $\C$  
such that $U\cap c(P)=\emptyset,$ $L\setminus \{t_0\}\subset U$
and let $p^{-1}_j(t),$ $1 \leq j \leq n,$ be the single-valued
branch of $P^{-1}(t)$ defined in $U$ by the
condition $p_j^{-1}\{l\setminus t_0\}=l_j\setminus\{a_j,b_{j}\}.$
If $k<n$ is a number such that $b_k=b$ then
for the same reason as above the function
$\phi(t)=\sum_{j=1}^k Q(p^{-1}_j(t))$ 
extends to an analytic in $U \cup D_L^+$ function but this fact does 
not imply now that $\phi(t)$ extends to an analytic in the whole
complex plane function since $D_L^+$ does not contain
$t_0\in c(P).$ Nevertheless,
if $V$ is a simply connected domain such that $U\subset V$ and
$t_0\notin V$ then $\phi(t)$ still extends to a single-valued analytic
function in $V.$ 
In particular, 
the analytic continuation of 
$\{U,\phi(t)\}$ along any simple closed curve $M$ such that 
$t_0\subset D_M^-$ coincides with
$\{U,\phi(t)\}.$

Let $t_1\in U$ be a point and let $M_1$ (resp. $M_2$) be a 
simple closed curve
such that $t_1\in M_1,$ $M_1\cap c(P)=\emptyset$
and $D^+_{M_1}\cap c(P)=t_0$ (resp. $t_1\in M_2,$ 
$M_2\cap c(P)=\emptyset$ and
$D^+_{M_2}\cap c(P)=c(P)\setminus t_0$). Define a permutation
$\rho_1\in S_n$ (resp. $\rho_2\in S_n$)
by the condition that the functional element
$\{U,p^{-1}_{\rho_1(j)}(t)\}$
(resp. $\{U,p^{-1}_{\rho_2(j)}(t)\}$) is the result
of the analytic continuation of the
functional element $\{U,p^{-1}_{j}(t)\},$ $1 \leq j \leq n,$ from $t_1$ 
along the curve $M_1$ (resp. $M_2$). 
Having in mind the identification of the set of elements
$\{U,p^{-1}_{j}(t)\},$
$1 \leq j \leq n,$ with the set of oriented edges of the graph
$\vec\Omega$ the permutations $\rho_1,\rho_2$ 
can be described as follows:
$\rho_1$ cyclically permutes the edges of $\vec\Omega$ 
around the vertices from which they go while cycles
$(j_1,j_2,...,j_k)$ of $\rho_2$ correspond to
simple cycles $(l_{j_1}, l_{j_2},..., l_{j_k})$
of the graph $\vec\Omega$ and 
$\rho_1\rho_2=(12...n)$ (see fig. 2).

To unload notation denote temporarily the element 
$\{U,Q(p^{-1}_{i}(t))\},$ $1\leq i \leq n,$ by $s_i.$
Since $t_0\subset D_{M_2}^-$ we have:
$$
0=\sum_{j=1}^{k} s_{\rho_2(j)}-\sum_{j=1}^{k} s_{j}=s_{\rho_2(k)}+\sum_{j=1}^{k-1}\left[
s_{\rho_2(j)}-
s_{j+1}\right]-s_{1}. \eqno(4)$$
Using $\rho_1\rho_2=(12...n)$ we can rewrite (4) as
$$s_{\rho^{-1}_1(k+1)}-s_{1}+\sum_{j=1}^{k-1}\left[
s_{\rho_2(j)}-s_{\rho_1\rho_2(j)}\right]=0.$$
Therefore, by the analytic continuation
$$s_{\rho_1^{f-1}(k+1)}-s_{\rho_1^f(1)}+\sum_{j=1}^{k-1}
\left[s_{\rho_1^{f}\rho_2(j)}-s_{\rho_1^{f+1}\rho_2(j)}
\right]=0 \eqno(5)$$  for $f\geq 0.$
Summing equalities (5) from $f=1$ to $f=o(\rho_1),$
where $o(\rho_1)$ is the order of the permutation $\rho_1,$
changing the order of summing,
and observing that 
$$\sum_{f=0}^{o(\rho_1)-1}\left[
s_{\rho_1^{f}\rho_2(j)}-
s_{\rho_1^{f+1}\rho_2(j)}\right]=
s_{\rho_2(j)}-s_{\rho_1^{o(\rho_1)}\rho_2(j)}=0$$ 
we obtain:
$$\sum_{f=0}^{o(\rho_1)-1}s_{\rho_1^f(k+1)}=
\sum_{f=0}^{o(\rho_1)-1}s_{\rho_1^f(1)}.\eqno(6)$$
If $a,b$ are not critical points
of $P(z)$ then
$p^{-1}_{\rho_1(1)}(t)=p^{-1}_{1}(t),$ 
$p^{-1}_{\rho_1(k+1)}(t)=p^{-1}_{k+1}(t)$
and (6) reduces to the equality
$Q(p^{-1}_{k+1}(t))=Q(p^{-1}_{1}(t)).$

Suppose now that at least one of the points $a,b$ is a critical 
point of $P(z).$
Observe that (6) is hold for any polynomial $Q(z)$ such that 
$q(z)=Q^{\prime}(z)$ satisfies
(*).Therefore, 
substituting in (6) $Q^j(z),$ $2\leq j \leq d_a+d_b-1,$ instead of $Q(z)$ we 
conclude that $$\sum_{s=0}^{o(\rho_1)-1}Q^j(p^{-1}_{\rho_1^s(k+1)}(t))=
\sum_{s=0}^{o(\rho_1)-1}Q^j(p^{-1}_{\rho_1^s(1)}(t))\eqno(7)$$
for all $j,$ $1\leq j \leq \d_b+d_b-1.$ Consider
a Vandermonde determinant $D=\parallel d_{j,i} \parallel,$
where $d_{j,i}=Q^j(p^{-1}_i(t)),$ 
$0\leq j \leq d_a+d_b-1$ and $i$ ranges the set of different indices
from the cycles of $\rho_1$ containing $1$ and $k+1.$
Since (7) implies that $D=0$ we conclude again that
$Q(p^{-1}_i(t))=Q(p^{-1}_j(t))$ for some $i\neq j,$ $1\leq i,j \leq n.$
Therefore, $P(z)$ and $Q(z)$ have a 
common right divisor in the composition algebra 
and we can finish the proof by the same argument as in  
section 2.3 taking into account that the multiplicity 
of a point $c\in \C$ with respect to $P(z)=\tilde P(W(z))$ is greater or
equal then the multiplicity of the point $W(c)$ with respect to 
$\tilde P(z).$

\pagebreak

\subsection{Proof of theorem 1.}
Suppose at first that $n=\deg P(z)$ is a prime number.
In this case the degree of the algebraic function 
$Q(P^{-1}(t))$ equals either $n$ or $1$
since $d(Q(P^{-1}(t)))$ divides $\deg P(z).$
If $d(Q(p^{-1}(t)))= n$ then 
Puiseux expansions at infinity
$$Q(p^{-1}_i(t))=\sum_{k\geq 
k_0}a_k\varepsilon^{ik}t^{\frac{k}{n}},\eqno(8)$$
$1\leq i \leq n,$ $a_k\in \C,$ $\varepsilon=exp(2\pi i/n),$
contain a coefficient $a_k\neq 0$ such that $k$ is not 
a multiple of $n.$
Substituting (8) in (6)
we conclude that $\varepsilon^k$ satisfies an algebraic 
polynomial with integer coefficients distinct from
the $n$-th cyclotomic polynomial $\Phi_n(z)=1+z+ ...
+z^{n-1}$.
Since $\varepsilon^k$ is a primitive $n$-th root of unity
it is a contradiction. Therefore, $d(Q(p^{-1}(t)))= 1$
and $Q(z)=\tilde Q(P(z))$ for some polynomial $\tilde Q(z)$.

Suppose now that $n$ is composite.
Since $P(z)$ is indecomposable
the group $G(P^{-1})$ is primitive by the Ritt theorem \cite{ri}.
By the Schur theorem (see e.g. \cite{wi}, Th. 25.3)
a primitive permutation group of composite degree $n$ 
which contains an $n$-cycle is doubly transitive. Therefore,
by lemma 2 equality (6) implies that
$Q(z)=\tilde Q(P(z))$ for some polynomial $\tilde Q(z)$. 

\vskip 0.2cm

\noindent{\bf Acknowledgments.} I am grateful to 
Y. Yomdin for drawing my attention to the polynomial moment
problem and for stimulating discussions.

\bibliographystyle{amsplain}

\end{document}